\documentclass{amsart}
\usepackage{amssymb}
\usepackage{mathrsfs}
\usepackage{graphicx}
\usepackage{comment}

\def\S{\mathcal S}

\newcommand{\nni}{\noindent}
\newcommand{\be}{\begin{equation}}
\newcommand{\ee}{\end{equation}}
\newcommand{\ba}{\begin{align}}
\newcommand{\ea}{\end{align}}

\newcommand{\abs}[1]{\lvert#1\rvert}

\newtheorem{example}{Example}[section]
\newtheorem{theorem}{Theorem}[section]
\newtheorem{corollary}[theorem]{Corollary}
\newtheorem{lemma}{Lemma}[section]
\newtheorem{alg}{Algorithm}[section]
\newtheorem{question}{Question}

\newenvironment{iolist}[1]%
{\begin{list}{}{%
\settowidth{\labelwidth}{\textsf{{\it #1.}}}%
\setlength{\labelsep}{4mm}%
\setlength{\leftmargin}{\labelwidth}%
\addtolength{\leftmargin}{\labelsep}%
}}%
{\end{list}}

\def\beq{\begin{equation}}\def\enq{\end{equation}}

\newenvironment{biglabellist}[1]%
{\begin{list}{}{%
\settowidth{\labelwidth}{\textsf{{\it #1.}}}%
\setlength{\labelsep}{2mm}%
\setlength{\leftmargin}{\labelwidth}%
\addtolength{\leftmargin}{\labelsep}%
\addtolength{\leftmargin}{4mm}%
\setlength{\itemsep}{6pt}%
\setlength{\listparindent}{0pt}%
\setlength{\topsep}{3pt}%
}}%

\title[group determinants]{The integer group determinants for $Q_{16}$}

\author[B. Paudel]{Bishnu Paudel}
\address{ Department of Mathematics\\
         Kansas State University\\
         Manhattan, KS 66506, USA}
\email{bpaudel@ksu.edu, pinner@math.ksu.edu}

\author[C. Pinner]{Christopher Pinner}

\date{\today}
\begin{document}

\begin{abstract}

We obtain a complete description of the integer group determinants for $Q_{16},$ the dicyclic or generalized quaternion group of order 16.

\end{abstract}

\maketitle

\section{Introduction}

At the meeting of the American Mathematical Society in Hayward, California, in April 1977, Olga Taussky-Todd \cite{TausskyTodd} asked whether one could characterize the values of the group determinant when the entries are all integers.
For a prime $p,$ a complete description was obtained for $\mathbb Z_{p}$ and $\mathbb Z_{2p}$, the cyclic groups of order $p$ and $2p$,  in \cite{Newman1} and \cite{Laquer}, and  for $D_{2p}$ and $D_{4p}$ the dihedral groups of order $2p$ and $4p$ in \cite{dihedral}. The values for $Q_{4n}$, the dicyclic group of order $4n$ were explored in \cite{dicyclic} 
with a near complete description for $Q_{4p}$.  In general though this quickly becomes a hard problem,
with only partial results known even for $\mathbb Z_{p^2}$ once $p\geq 7$ (see \cite{Newman2} and \cite{Mike}).

The remaining groups of order less than 15 were tackled in  \cite{smallgps} and $\mathbb Z_{15}$ in \cite{bishnu1}. 
The integer group determinants have been determined for  all  five  abelian groups of order 16 ($\mathbb Z_2 \times \mathbb Z_8$, $\mathbb Z_{16}$, $\mathbb Z_2^4$, $\mathbb Z_4^2$, $\mathbb Z_2^2 \times\mathbb Z_4$  in  \cite{Yamaguchi1,Yamaguchi2,Yamaguchi3,Yamaguchi4,Yamaguchi5}), and for three of the  non-abelian groups
($D_{16}$, $\mathbb Z_2\times D_8$, $\mathbb Z_2 \times Q_8$ in \cite{dihedral,ZnxH}).

Here we determine the the group determinants for $Q_{16}$,  the dicyclic or generalized quaternion group of order 16.
$$ Q_{16}=\langle X,Y \; | \; X^8=1,\; Y^2=X^4,\; XY=YX^{-1}\rangle. $$
 This leaves five  unresolved non-abelian groups of order 16

\begin{theorem} The even integer group determinants for $Q_{16}$ are exactly the multiples of $2^{10}$.

The odd integer group determinants are all the integers $n\equiv 1$ mod 8 plus  those $n\equiv 5$ mod 8 of the form
$n=mp^2$ where $m\equiv 5$ mod 8 and $p\equiv 7$ mod $8$ is prime.

\end{theorem}

We shall think here of the  group determinant as being defined on elements of the group ring $\mathbb Z [G]$
$$ \mathcal{D}_G\left( \sum_{g\in G} a_g g \right)=\det\left( a_{gh^{-1}}\right) .$$
Frobenius \cite{Frob} observed that the group determinant can be factored using the  groups representations (see for example \cite{Conrad} or \cite{book})
and an explicit  expression for a dicyclic group determinant was given in \cite{smallgps}. For $Q_{16}$, arranging the
16 coefficients into two polynomials of degree 7
$$ f(x)=\sum_{j=0}^7 a_j x^j,\;\; g(x)=\sum_{j=0}^7 b_jx^j, $$
and writing the primitive 8th root of unity $\omega:=e^{2\pi i/8}=\frac{\sqrt{2}}{2}(1+i)$, this becomes
\be \label{form}\mathcal{D}_G\left( \sum_{j=0}^7 a_j X^j + \sum_{j=0}^7 b_j YX^j\right)  =ABC^2D^2 \ee
with integers $A,B,C,D$ from
\begin{align*}
A=& f(1)^2- g(1)^2\\
B=& f(-1)^2-g(-1)^2\\
C=& |f(i)|^2-|g(i)|^2 \\
D=& \left(|f(\omega)|^2+|g(\omega)|^2\right)\left(|f(\omega^3)|^2+|g(\omega^3)|^2\right).
\end{align*}
From \cite[Lemma 5.2]{dicyclic} we know that the even values must be multiples of $2^{10}$.  The odd values must be
1 mod 4 (plainly $f(1)$ and $g(1)$ must be of opposite parity and $A\equiv B\equiv \pm 1$  mod 4 with $(CD)^2\equiv 1$ mod 4).

\section{Achieving the values $n\not \equiv 5$ mod 8}

We can achieve all the multiples of $2^{10}$. 

Writing $h(x):=(x+1)(x^2+1)(x^4+1),$ we achieve the $2^{10}(-3+4m)$ from
$$
f(x) = (1-m)h(x),\quad
g(x)=1+x^2+x^3+x^4-mh(x), $$
the $2^{10}(-1+4m)$  from
$$ f(x)= 1+x+x^4+x^5-mh(x),\;\;\;\;
g(x)= 1+x-x^3-x^7-mh(x), $$
the $2^{11}(-1+2m)$ from
$$ f(x)= 1+x+x^2+x^3+x^4+x^5-mh(x),\;\;\quad
g(x)=1+x^4-mh(x), $$
and the $2^{12}m$ from 
$$ f(x)= 1+x+x^4+x^5-x^6-x^7-mh(x),\;\;
g(x)= 1+x-x^3+x^4+x^5-x^7+mh(x). $$

We can achieve all the $n\equiv  1$ mod 8; the $1+16m$ from
$$ f(1)=1+mh(x),\;\; g(x)=mh(x), $$
and the $-7+16m$ from 
$$f(x)= 1-x+x^2+x^3+x^7- mh(x),\;\;
g(x)= 1+x^3+x^4+x^7-mh(x). $$

\section{ The form of the $n\equiv 5$ mod 8}

This leaves the $n\equiv 5$ mod 8. Since $(CD)^2\equiv 1$ mod 8 we must have $AB\equiv 5$ mod 8. Switching $f$ and $g$ as necessary we assume that $f(1),f(-1)$ are odd and $g(1),g(-1)$ even. Replacing $x$ by $-x$  if needed we can assume that $g(1)^2\equiv 4$ mod 8 and $g(-1)^2\equiv 0$ mod 8.

We write 
$$ F(x)=f(x)f(x^{-1})= \sum_{j=0}^7 c_j (x+x^{-1})^j, \quad  G(x)=g(x)g(x^{-1})= \sum_{j=0}^7 d_j (x+x^{-1})^j, $$
with the $c_j,d_j$ in $\mathbb Z$.

From $F(1),F(-1)\equiv 1$ mod 8 we have
$$ c_0+2c_1+4c_2 \equiv 1 \text{ mod }8, \quad  c_0-2c_1+4c_2 \equiv 1 \text{ mod }8,  $$
and $c_0$ is odd and $c_1$ even.
From $G(1)\equiv 4$, $G(-1)\equiv 0$ mod 8 we have
$$  d_0+2d_1+4d_2 \equiv 4 \text{ mod 8}, \quad d_0-2d_1+4d_2 \equiv 0 \text{ mod } 8, $$
and $d_0$ is even and $d_1$ is odd. 
Since $\omega+\omega^{-1}=\sqrt{2}$ we get
\begin{align*} F(\omega)  & = (c_0+2c_2+4c_4+\ldots ) + \sqrt{2}(c_1+2c_3+4c_5+\cdots),\\
G(\omega)  & = (d_0+2d_2+4d_4+\ldots ) + \sqrt{2}(d_1+2d_3+4d_5+\cdots),
\end{align*}
and
$$|f(\omega)|^2+|g(\omega)|^2=  F(\omega)+G(\omega) = X+ \sqrt{2} Y>0, \;\; \quad X, Y \text{odd, }  $$
with $ |f(\omega^3)|^2+|g(\omega^3)|^2=F(\omega^3)+G(\omega^3) = X- \sqrt{2} Y>0$. Hence  the positive integer $D=X^2-2Y^2\equiv -1$ mod 8.
Notice that primes 3 and 5 mod 8 do not split in $\mathbb Z[\sqrt{2}]$ so only their squares can occur in $D$. Hence
$D$ must contain at least one prime $p\equiv 7$ mod 8, giving the claimed form of the values 5 mod 8.

\section{Achieving the specified values 5 mod 8}
Suppose that $p\equiv 7$ mod 8 and $m\equiv 5$ mod 8. We need to achieve $mp^2$.

Since $p\equiv 7$ mod 8 we know that $\left(\frac{2}{p}\right)=1$ and $p$ splits in $\mathbb Z[\sqrt{2}].$ Since $\mathbb Z[\sqrt{2}]$  is 
a UFD, a generator for the prime factor gives a  solution to
$$ X^2-2Y^2=p, \;\; X,Y\in \mathbb N. $$
Plainly $X,Y$ must both be odd and $X+\sqrt{2}Y$ and $X-\sqrt{2}Y$ both positive. 
Since $(X+\sqrt{2}Y)(3+2\sqrt{2})=(3X+4Y)+\sqrt{2}(2X+3Y)$ there will be $X,Y$ with $X\equiv 1$ mod 4 and with
$X\equiv -1$ mod 4.

Cohn \cite{Cohn} showed that $a+b\sqrt{2}$ in $\mathbb Z[\sqrt{2}]$ is a sum of four squares in $\mathbb Z[\sqrt{2}]$ if and only if $2\mid b$. Hence we can write
$$ 2(X+\sqrt{2}Y)= \sum_{j=1}^4 (\alpha_j + \beta_j\sqrt{2})^2, \;\;\alpha_j,\beta_j\in \mathbb Z. $$
That is,
$$ 2X=\sum_{j=1}^4 \alpha_j^2+ 2\sum_{j=0}^4 \beta_j^2,\;\;\quad Y=\sum_{j=1}^4\alpha_j\beta_j.$$
Since $Y$ is odd we must have at least one pair, $\alpha_1$, $\beta_1$ say,  both odd. Since $2X$ is even we must have two
or four of the $\alpha_i$ odd. Suppose that $\alpha_1$, $\alpha_2$ are odd and $\alpha_3,\alpha_4$ have the same parity.
We get
\begin{align*} X+\sqrt{2}Y & = \left(  \frac{\alpha_1+\alpha_2}{2} + \frac{\sqrt{2}}{2}(\beta_1+\beta_2)\right)^2+ \left(  \frac{\alpha_1-\alpha_2}{2} + \frac{\sqrt{2}}{2}(\beta_1-\beta_2)\right)^2 \\
  & \quad +  \left(  \frac{\alpha_3+\alpha_4}{2} + \frac{\sqrt{2}}{2}(\beta_3+\beta_4)\right)^2+ \left(  \frac{\alpha_3-\alpha_4}{2} + \frac{\sqrt{2}}{2}(\beta_3-\beta_4)\right)^2.
\end{align*} 
Writing 
$$ f(\omega)=a_0+a_1\omega+a_2\omega^2+a_3\omega^3=a_0+ \frac{\sqrt{2}}{2}(1+i)a_1+a_2i+ \frac{\sqrt{2}}{2}(-1+i)a_3,$$
we have
$$ \abs{f(\omega)}^2 =\left(a_0+ \frac{\sqrt{2}}{2}(a_1-a_3)\right)^2  +  \left(a_2+ \frac{\sqrt{2}}{2}(a_1+a_3)\right)^2  $$
and can make
$$ |f(\omega)|^2+|g(\omega)|^2 = X + \sqrt{2}Y $$
with the selection of integer coefficients for $f(x)=\sum_{j=0}^3a_jx^j$ and $g(x)=\sum_{j=0}^3 b_jx^j$
\begin{align*}  a_0=&\frac{1}{2}(\alpha_1-\alpha_2),\quad  a_1 =\beta_1,\quad   a_2=\frac{1}{2}(\alpha_1+\alpha_2), \quad a_3=\beta_2, \\
    b_0=& \frac{1}{2}(\alpha_3-\alpha_4),\quad  b_1 =\beta_3,\quad   b_2=\frac{1}{2}(\alpha_3+\alpha_4), \quad b_3=\beta_4. 
\end{align*}
These $f(x)$, $g(x)$ will then give $D=p$ in \eqref{form}.
We can also determine the parity of the coefficients. 

\vskip0.1in
\nni
{\bf Case 1}: the $\alpha_i$ are all odd.

Notice that $a_0$ and $a_2$ have opposite parity, as do $b_0$ and $b_2$. Since $Y$ is odd we must have one or three of the 
$\beta_i$ odd.

 If $\beta_1$ is odd and $\beta_2,\beta_3,\beta_4$ all even, then $2X\equiv 6$ mod 8 and $X\equiv -1$ mod 4.
Then $a_0,a_1,a_2,a_3$ are either odd, odd, even, even or even, odd, odd, even and $f(x)=u(x)+2k(x)$
with $u(x)=1+x$ or $x(1+x)$. Likewise $b_0,b_1,b_2,b_3$  are odd, even, even, even or even, even, odd, even
and $g(x)=v(x)+2s(x)$ with $v(x)=1$ or $x^2$. Hence if we take
\be \label{shift}  f(x)=u(x)+(1-x^4)k(x)-mh(x),\quad g(x)=v(x)+(1-x^4)s(x)-mh(x), \ee
we get $A=3-16m$, $B=-1$, $C=1$, $D=p$ and we achieve $(16m-3)p^2$ in \eqref{form}.

If three $\beta_i$ are odd then $2X\equiv 10$ mod 8 and $X\equiv 1$ mod 4. We assume $\beta_1,\beta_2,\beta_3$ are 
odd and $\beta_4$ even. Hence  $a_0,a_1,a_2,a_3$ are either odd, odd, even, odd or even, odd, odd, odd and 
$f(x)=u(x)+2k(x)$ with $u(x)=1+x+x^3$ or $x(1+x+x^2)$ and $b_0,b_1,b_2,b_3$  are odd, odd, even, even or even, odd, odd, even
and $g(x)=v(x)+2s(x)$ with $v(x)=1+x$ or $x(1+x)$. In this case \eqref{shift} gives
$A=(5-16m)$, $B=1$, $C=-1$, $D=p$ achieving $(5-16m)p^2$.

\vskip0.1in
\nni
{\bf Case 2}: $\alpha_1$, $\alpha_2$ are odd, $\alpha_3$, $\alpha_4$ are even.

In this case $a_0$, $a_2$ will have opposite parity and $b_0$, $b_2$ the same parity.
Since $Y$ is odd we must have $\beta_1$ odd, $\beta_2$ even. Since $2X\equiv 2$ mod 4 we must have one more odd $\beta_i$, say $\beta_3$ odd and $\beta_4$ even. 
If $\alpha_3\equiv \alpha_4$ mod 4 then $2X\equiv 6$ mod 8 and $X\equiv -1$ mod 4. Hence
 $a_0,a_1,a_2,a_3$ are either odd, odd, even, even or even, odd, odd, even, that is $u(x)=1+x$ or $x(1+x)$  and  $b_0,b_1,b_2,b_3$  are even, odd, even, even and $v(x)=x^2$ and again \eqref{shift} gives $(16m-3)p^2$.

If $\alpha_3\not\equiv \alpha_4$ mod 4 then $2X\equiv 10$ mod 8 and $X\equiv 1$ mod 4. In this case
 $a_0,a_1,a_2,a_3$ are either odd, odd, even, even or even, odd, odd, even, that is $u(x)=1+x$ or $x(1+x)$  and  $b_0,b_1,b_2,b_3$  are odd, odd, odd, even and $v(x)=1+x+x^2$ and again \eqref{shift} gives $(5-16m)p^2$.

Hence, in either case, starting with an $X\equiv 1$ mod 4 gives the $mp^2$ with $m\equiv 5$ mod 16 and 
an  $X\equiv -1$ mod 4 the $mp^2$ with $m\equiv -3$ mod 16.

\section*{Acknowledgement} 
\nni
We thank Craig Spencer for directing us to Cohn's four squares theorem in $\mathbb Z[\sqrt{2}]$.


\begin{thebibliography}{99}
\bibitem{dihedral}
T. Boerkoel \& C. Pinner, \textit{Minimal  group determinants and the  Lind-Lehmer problem for dihedral groups},   Acta Arith. \textbf{186} (2018), no. 4, 377-395.	arXiv:1802.07336 [math.NT].

\bibitem{Cohn} 
H.\ Cohn, \textit{Decomposition into four integral squares in the fields of $2^{\frac{1}{2}}$ and $3^{\frac{1}{2}}$,} 
Amer. J. Math. \textbf{82} (1960), 301-322.

\bibitem{Conrad}
K.\ Conrad, \textit{The origin of representation theory},  Enseign. Math. (2) \textbf{44} (1998), no. 3-4, 361-392.


\bibitem{Frob} F.\ G.\ Frobenius, \textit{\"{U}ber die  Primefactoren  der  Gruppendeterminante},  Gesammelte Ahhand-lungen, Band III, Springer, New York, 1968, pp. 38–77. MR0235974


\bibitem{book}
K.\ Johnson, Group Matrices Group Determinants and Representation Theory, Lecture Notes in Mathematics 2233,  Springer 2019.


\bibitem{Mike}
M.\ Mossinghoff and C.\ Pinner, \textit{Prime power order circulant determinants}, (arxiv 2205.12439v2.)


\bibitem{Newman1}
M.\ Newman, \textit{ On a problem suggested by Olga Taussky-Todd}, Ill. J. Math. \textbf{24} (1980), 156-158.

\bibitem{Newman2}
M.\ Newman, \textit{Determinants of circulants of prime power order}, Linear Multilinear Algebra \textbf{ 9}
(1980), no. 3, 187–191. MR0601702.




\bibitem{dicyclic}
B.\ Paudel and C.\ Pinner, \textit{Minimal group determinants for dicyclic groups,} Mosc. J. Comb. Number Theory \textbf{10} (2021), no.3, 235-248.

\bibitem{bishnu1}
B.\ Paudel and C.\ Pinner, \textit{Integer circulant determinants of order 15,} Integers \textbf{22} (2022), Paper No. A4.

\bibitem{ZnxH}
B.\ Paudel and C.\ Pinner, \textit{The group determinants for $\mathbb Z_n x H$,}  arXiv:2211.09930 [math.NT].

\bibitem{Laquer}
H. T. Laquer, \textit{ Values of circulants with integer entries,}  in A Collection of Manuscripts Related
to the Fibonacci Sequence, Fibonacci Assoc., Santa Clara, 1980, pp. 212–217. MR0624127.

\bibitem{smallgps}
C. Pinner and C. Smyth, \textit{Integer group determinants for small groups},  Ramanujan J. \textbf{51} (2020), no. 2, 421-453.

\bibitem{TausskyTodd}
O. Taussky Todd, \textit{Integral group matrices}, Notices Amer. Math. Soc. 24 (1977), no. 3, A-345.
Abstract no. 746-A15, 746th Meeting, Hayward, CA, Apr. 22–23, 1977.


\bibitem{Yamaguchi1}
Y.\ Yamaguchi and N.\ Yamaguchi, 
\textit{Generalized Dedekind’s theorem and its application to integer group determinants}, 2022. arXiv:2203.14420v2 [math.RT].


 
\bibitem{Yamaguchi2}
Y.\ Yamaguchi and N.\ Yamaguchi, 
\textit{Integer circulant determinants of order 16},  Ramanujan J., arXiv:2204.05014 [math.NT].

\bibitem{Yamaguchi3}
Y.\ Yamaguchi and N.\ Yamaguchi, 
\textit{Integer group determinants for $C_2^4$}, 2022.
arXiv:2203.14420v2 [math.RT]


\bibitem{Yamaguchi4}
Y.\ Yamaguchi and N.\ Yamaguchi, \textit{Integer group determinants for $C_4^2$}, arXiv:2211.01597 [math.NT].

\bibitem{Yamaguchi5}
Y.\ Yamaguchi and N.\ Yamaguchi, \textit{Integer group determinants for abelian groups of order 16}, arXiv:2211.14761 [math.NT].
\end{thebibliography}
\end{document}